\documentclass[11pt,a4paper]{article}

\usepackage[margin=1in]{geometry}
\usepackage{amsmath,amssymb,amsthm,mathtools}
\usepackage{mathrsfs}
\usepackage{enumitem}
\usepackage{microtype}
\usepackage{color}
\usepackage[colorlinks=true,linkcolor=blue,citecolor=blue,urlcolor=blue]{hyperref}

\theoremstyle{plain}
\newtheorem{theorem}{Theorem}[section]
\newtheorem{lemma}[theorem]{Lemma}
\newtheorem{corollary}[theorem]{Corollary}

\newtheorem{conj}{Conjecture}[section]
\theoremstyle{definition}

\theoremstyle{remark}
\newtheorem{remark}[theorem]{Remark}

\newcommand{\qrad}[1]{q(#1)}

\title{Maximum signless Laplacian spectral radius of graphs with given minimum degree and edge connectivity}
\author{Joyentanuj Das\thanks{Department of Mathematics, SRM Institute of Science and Technology Chennai, Tamil Nadu, India. Email: joyentanuj@gmail.com, joyentad@srmist.edu.in.}\  \and Iswar Mahato \thanks{Department of Mathematics, SRM University-AP, Andhra Pradesh 522240, India. Email: iswarmahato02@gmail.com, iswar.m@srmap.edu.in.}}
\date{}

\begin{document}

\maketitle

\begin{abstract}
Let $G$ be a simple connected graph of order $n$ with minimum degree $\delta$ and edge connectivity $k$. Recently, in [Linear Algebra Appl. 751 (2026) 318-330], Zhao, Ye, and Fan characterized the extremal graph with maximum signless Laplacian spectral radius among all graphs of order $n$ with minimum degree  $\delta$ and edge connectivity $k$ when $n \geq \delta^2+4\delta+3$ and $k \in \{0,1,2\}$, and proposed the corresponding problem for $3 \leq k \leq \delta-1$. In this paper, we solve the problem for the entire range $0 \leq k \leq \delta-1$ and completely remove the order restriction. Specifically, we prove that 
\[
\qrad{G} \leq \qrad{B_{n,n-\delta-1}^{k}},
\]
with equality if and only if $G \cong B_{n,n-\delta-1}^{k}$, where $B_{n,n-\delta-1}^{k}$ is the graph obtained from $K_{n-\delta-1} \cup K_{\delta+1}$ by adding $k$ edges joining a single vertex of $K_{n-\delta-1}$ to $k$ distinct vertices of $K_{\delta+1}$. 
\end{abstract}


\section{Introduction and Preliminaries}

\subsection{Introduction}
Throughout this paper, all graphs are finite, simple and undirected. For a graph $G$, let $A(G)$ and $D(G)$ denote its adjacency matrix and diagonal degree matrix, respectively. The signless Laplacian matrix of $G$ is defined as
\[
Q(G)=D(G)+A(G),
\]
and its largest eigenvalue, denoted by $\qrad{G}$, is called the signless Laplacian spectral radius or the $Q$-index of $G$. The spectral properties of the signless Laplacian have been extensively studied; see, for instance, the foundational works \cite{CvetkovicRowlinsonSimic2007, Das} and the comprehensive monographs \cite{BrouwerHaemers, CvetkovicRowlinsonSimic, CioabaGregoryHaemers}.

The edge connectivity $\kappa'(G)$ of a connected graph $G$ is the minimum number of edges whose deletion disconnects $G$. The fundamental relation $\kappa'(G)\leq\delta(G)$, where $\delta(G)$ is the minimum degree, makes it natural to ask how the $Q$-index behaves when the minimum degree and edge connectivity are prescribed. This belongs to a broad class of extremal spectral problems in which a structural parameter is controlled by the largest eigenvalue of a graph matrix. For related results on the adjacency spectral radius and edge connectivity, we refer to \cite{Guo}.

For positive integers $n,b,k$ with $k\leq n-b$, let $B_{n,b}^{k}$ be the graph obtained from $K_b\cup K_{n-b}$ by choosing a vertex $v\in V(K_b)$ and adding $k$ edges from $v$ to $k$ distinct vertices of $K_{n-b}$. Recently, Zhao, Ye and Fan \cite{ZhaoYeFan} studied the maximum $Q$-index of graphs with prescribed minimum degree and edge connectivity. Under the restrictive hypothesis $n\geq\delta^2+4\delta+3$, they reduced the problem for $0\leq k\leq\delta-1$ to graphs formed from $K_{\delta+1}$ and $K_{n-\delta-1}$ by inserting $k$ joining edges. They settled the cases $k\in\{0,1,2\}$ and proposed the following conjecture for $3\leq k\leq\delta-1$:

\begin{conj}
  For $3\leq k \leq \delta-1$, the graph $B_{n,n-\delta-1}^{k}$ is the extremal graph with maximum signless Laplacian spectral radius among all graphs of order $n$ with minimum degree $\delta$ and edge connectivity $k$. 
\end{conj}

In this paper, we answer this question affirmatively and completely remove the order restriction. Furthermore, we derive the precise asymptotic expansion of the maximum $Q$-index as $n \to \infty$, revealing a striking hierarchy in how the structural parameters influence the spectral radius at large orders.

Our proof of the extremal result has two principal ingredients. First, among all graphs obtained by inserting $k$ edges between two fixed cliques, the $Q$-index is maximized when the joining edges form a star centered in the larger clique. Second, transferring a vertex from the smaller clique to the larger clique strictly increases the $Q$-index of this extremal graph. These facts, combined with a minimum-edge-cut argument, yield the complete structural result. The asymptotic analysis is then carried out by expanding the exact characteristic equation of the extremal graph.

\subsection{Preliminaries}
For disjoint sets $S, T \subseteq V(G)$, let $E_G(S,T)$ be the set of edges with one endpoint in $S$ and the other in $T$, and write $e_G(S,T) = |E_G(S,T)|$. We will frequently use the following standard facts.

\begin{lemma}\label{lem:monotonicity}
Let $G$ be a connected graph and let $H$ be a spanning subgraph of $G$. Then $\qrad{H} \leq \qrad{G}$. If $H$ is connected, equality holds if and only if $H = G$.
\end{lemma}
\begin{proof}
Every edge $uv$ contributes the positive semidefinite matrix $(\mathbf{e}_u + \mathbf{e}_v)(\mathbf{e}_u + \mathbf{e}_v)^\top$ to the signless Laplacian. The assertion follows directly from the Rayleigh principle and, in the equality case, the strict positivity of the Perron vector.
\end{proof}

\begin{lemma}\label{lem:cut-size}
Let $G$ have minimum degree $\delta$. If $\emptyset \neq S \subsetneq V(G)$ and $e_G(S, V(G) \setminus S) \leq \delta-1$, then $|S| \geq \delta+1$.
\end{lemma}
\begin{proof}
Suppose for a contradiction that $|S| \leq \delta$. Then every vertex $u \in S$ has at most $|S|-1$ neighbors inside $S$, meaning it has at least $\delta - |S| + 1$ neighbors outside $S$. Consequently,
\[
e_G(S, V(G) \setminus S) \geq |S|(\delta - |S| + 1).
\]
The function $f(x) = x(\delta - x + 1)$ is minimized on the interval $[1, \delta]$ at the boundaries. For $x=1$, $f(1) = \delta$. For $x=\delta$, $f(\delta) = \delta$. Thus, $e_G(S, V(G) \setminus S) \geq \delta$, which contradicts the assumption.
\end{proof}

We also require the resolvent of the signless Laplacian of a complete graph. For $m \geq 1$, $Q(K_m) = (m-2)I_m + J_m$, where $J_m$ is the all-ones matrix. Thus, whenever $x > 2m-2$, the inverse $(xI_m - Q(K_m))^{-1}$ can be explicitly computed using the Sherman-Morrison formula:
\begin{equation}\label{eq:resolvent-clique}
\bigl(xI_m - Q(K_m)\bigr)^{-1} = \alpha_m(x)I_m + \gamma_m(x)J_m,
\end{equation}
where
\begin{equation}\label{eq:alpha-gamma}
\alpha_m(x) = \frac{1}{x - m + 2}, \qquad \gamma_m(x) = \frac{1}{(x - m + 2)(x - 2m + 2)}.
\end{equation}

\section{Main results}

In this section, we establish the structural lemmas that lead to the main extremal theorem. We first determine the optimal arrangement of a prescribed number of edges between two fixed cliques.

\begin{lemma}\label{lem:concentration}
Let $b \geq a \geq k \geq 1$, and let $H$ be obtained from $K_a \cup K_b$ by adding $k$ edges between the two cliques. Then
\[
\qrad{H} \leq \qrad{B_{a+b,b}^{k}}.
\]
If $b > a$, equality holds if and only if the $k$ joining edges have a common endpoint in $K_b$ and distinct endpoints in $K_a$. If $a = b$, equality holds if and only if the joining edges form a star in one of the two cliques.
\end{lemma}
\begin{proof}
Write the joining edges as $e_i = u_i v_i$, where $u_i \in V(K_a)$ and $v_i \in V(K_b)$ for $1 \leq i \leq k$. Set $Q_0 = Q(K_a \cup K_b)$ and $z_i = \mathbf{e}_{u_i} + \mathbf{e}_{v_i}$. Let $Z$ be the $(a+b) \times k$ matrix with columns $z_1, \ldots, z_k$. Then
\begin{equation}\label{eq:rank-k}
Q(H) = Q_0 + Z Z^\top.
\end{equation}
Since $H$ is a connected proper supergraph of $K_b$, we have $\qrad{H} > 2b-2$. For $x > 2b-2$, put $R(x) = (xI - Q_0)^{-1}$. Define $k \times k$ matrices $P_a$ and $P_b$ by $(P_a)_{ij} = \mathbf{1}_{\{u_i = u_j\}}$ and $(P_b)_{ij} = \mathbf{1}_{\{v_i = v_j\}}$. Using \eqref{eq:resolvent-clique}, a direct calculation yields the compressed resolvent:
\begin{equation}\label{eq:compressed-resolvent}
Z^\top R(x) Z = \bigl(\gamma_a(x) + \gamma_b(x)\bigr)J_k + \alpha_a(x)P_a + \alpha_b(x)P_b.
\end{equation}
Since $b \geq a$, we have $\alpha_b(x) \geq \alpha_a(x) > 0$. Let $y = (y_1, \ldots, y_k)^\top \geq 0$. Notice that
\begin{align*}
y^\top (P_a - I)y &= 2 \sum_{\substack{i < j \\ u_i = u_j}} y_i y_j, \\
y^\top (J - P_b)y &= 2 \sum_{\substack{i < j \\ v_i \neq v_j}} y_i y_j.
\end{align*}
Because the joining edges are distinct, $u_i = u_j$ implies $v_i \neq v_j$. It follows that for any $y \geq 0$,
\begin{align}
&y^\top \bigl(\alpha_a I + \alpha_b J - \alpha_a P_a - \alpha_b P_b\bigr)y \notag \\
&\quad = \alpha_b y^\top (J - P_b)y - \alpha_a y^\top (P_a - I)y \geq 0. \label{eq:key-quadratic}
\end{align}
Let $H^\star = B_{a+b,b}^{k}$. For its joining edges, $P_a^\star = I$ and $P_b^\star = J$. Taking a nonnegative Perron vector in \eqref{eq:key-quadratic} gives
\begin{equation}\label{eq:resolvent-rho-comparison}
\rho\bigl(Z^\top R(x) Z\bigr) \leq \rho\bigl((Z^\star)^\top R(x) Z^\star\bigr).
\end{equation}
Let $q = \qrad{H}$ and let $X > 0$ be a Perron vector of $Q(H)$. From \eqref{eq:rank-k}, $(qI - Q_0)X = Z Z^\top X$. With $y = Z^\top X > 0$, this becomes $y = Z^\top R(q) Z y$. The matrix $Z^\top R(q) Z$ is positive, so its spectral radius equals one. Consequently, \eqref{eq:resolvent-rho-comparison} implies $\rho\bigl((Z^\star)^\top R(q) Z^\star\bigr) \geq 1$. The expression on the left is strictly decreasing in $q$ and equals one at $q = \qrad{H^\star}$. Therefore, $\qrad{H} \leq \qrad{H^\star}$.

Suppose first that $b > a$. Then $\alpha_b(x) > \alpha_a(x)$. Since the Perron vector is strictly positive, equality in \eqref{eq:key-quadratic} is possible only when no two joining edges share their endpoint in $K_a$ (so $P_a = I$) and every two joining edges share their endpoint in $K_b$ (so $P_b = J$). Hence, the joining edges form the asserted star. When $a = b$, equality requires every pair of joining edges to have a common endpoint in at least one of the cliques. A pairwise intersecting family of edges in a bipartite graph is a star, completing the proof.
\end{proof}

The next lemma shows that making the two clique orders more unequal strictly increases the extremal $Q$-index.

\begin{lemma}\label{lem:transfer}
Let $b \geq a$ and $1 \leq k \leq a-2$. Then
\[
\qrad{B_{a+b,b}^{k}} < \qrad{B_{a+b,b+1}^{k}},
\]
where the graph on the right has clique orders $a-1$ and $b+1$.
\end{lemma}
\begin{proof}
For $x > 2b-2$, define $\Phi_{a,b}(x) = \alpha_a(x) + k\gamma_a(x) + k\alpha_b(x) + k\gamma_b(x)$. For $B_{a+b,b}^{k}$, the compressed resolvent in \eqref{eq:compressed-resolvent} evaluates to $\alpha_a(x)I_k + \bigl(\alpha_b(x) + \gamma_a(x) + \gamma_b(x)\bigr)J_k$. It follows that $q_{a,b} = \qrad{B_{a+b,b}^{k}}$ is the unique solution, in $x > 2b-2$, of $\Phi_{a,b}(x) = 1$.

If $q_{a,b} \leq 2b$, then $K_{b+1}$ is a proper subgraph of $B_{a+b,b+1}^{k}$, and hence $\qrad{B_{a+b,b+1}^{k}} > 2b \geq q_{a,b}$. It remains to consider $q_{a,b} > 2b$.

For fixed $x > 2b$, regard $\alpha_m(x)$ and $\gamma_m(x)$ as functions of the continuous real variable $m$. They are strictly increasing and strictly convex in the range under consideration. Indeed,
\[
\frac{\partial^2 \alpha_m(x)}{\partial m^2} = \frac{2}{(x + 2 - m)^3} > 0.
\]
Writing $c = x + 2$, we have $\gamma_m(x) = \frac{2}{c(c-2m)} - \frac{1}{c(c-m)}$, and therefore
\[
\frac{\partial^2 \gamma_m(x)}{\partial m^2} = \frac{16}{c(c-2m)^3} - \frac{2}{c(c-m)^3} > 0.
\]
At $x = q_{a,b}$, set $\Delta \alpha_m = \alpha_{m+1} - \alpha_m$ and $\Delta \gamma_m = \gamma_{m+1} - \gamma_m$. Strict convexity and $b > a-1$ imply $\Delta \alpha_b > \Delta \alpha_{a-1}$ and $\Delta \gamma_b > \Delta \gamma_{a-1}$. Consequently,
\begin{align*}
\Phi_{a-1,b+1}(q_{a,b}) - \Phi_{a,b}(q_{a,b}) &= k(\Delta \alpha_b - \Delta \alpha_{a-1}) + k(\Delta \gamma_b - \Delta \gamma_{a-1}) > 0.
\end{align*}
Since $\Phi_{a,b}(q_{a,b}) = 1$, we have $\Phi_{a-1,b+1}(q_{a,b}) > 1$. Because $\Phi_{a-1,b+1}$ is strictly decreasing on $(2b, \infty)$, its unique root must be larger than $q_{a,b}$. This proves the claim.
\end{proof}

We can now solve the proposed problem in its full range and without any order restriction.

\begin{theorem}\label{thm:main}
Let $G$ be a graph of order $n$, minimum degree $\delta$, and edge connectivity $1 \leq k \leq \delta-1$. Then
\begin{equation}\label{eq:main-bound}
\qrad{G} \leq \qrad{B_{n,n-\delta-1}^{k}}.
\end{equation}
Equality holds if and only if $G \cong B_{n,n-\delta-1}^{k}$.
\end{theorem}
\begin{proof}
Choose $S \subset V(G)$ such that $e_G(S, V(G) \setminus S) = k$. Put $T = V(G) \setminus S$, $s = |S|$, and $t = |T|$, and assume without loss of generality that $s \leq t$. Since $k \leq \delta-1$, Lemma~\ref{lem:cut-size}, applied to both sides of the cut, gives $s \geq \delta+1$ and $t \geq \delta+1$. Thus, the existence of $G$ inherently forces $n \geq 2\delta+2$.

Let $H$ be obtained from $G$ by adding all missing edges inside $S$ and inside $T$. Then $H[S] \cong K_s$, $H[T] \cong K_t$, and $H$ has precisely the original $k$ edges between $S$ and $T$. By Lemmas~\ref{lem:monotonicity} and \ref{lem:concentration},
\begin{equation}\label{eq:chain-one}
\qrad{G} \leq \qrad{H} \leq \qrad{B_{n,t}^{k}}.
\end{equation}
If $s > \delta+1$, then $s \geq \delta+2$, which implies $k \leq \delta-1 \leq s-3 < s-2$. Lemma~\ref{lem:transfer} can therefore be applied repeatedly, giving
\begin{equation}\label{eq:repeated-transfer}
\qrad{B_{n,t}^{k}} < \qrad{B_{n,t+1}^{k}} < \cdots < \qrad{B_{n,n-\delta-1}^{k}}.
\end{equation}
Together with \eqref{eq:chain-one}, this proves that equality in \eqref{eq:main-bound} is impossible when $s > \delta+1$. Hence, equality requires $s = \delta+1$ and $t = n-\delta-1$.

Equality in the first inequality of \eqref{eq:chain-one} forces $G = H$, meaning $G[S]$ and $G[T]$ must already be complete graphs. Equality in the second inequality forces the $k$ cut edges to form a star centered in $K_t$ with distinct endpoints in $K_s$. Thus, $G \cong B_{n,n-\delta-1}^{k}$. Conversely, $B_{n,n-\delta-1}^{k}$ has order $n$, minimum degree $\delta$, and edge connectivity $k$, and therefore attains the bound.
\end{proof}

For the disconnected case ($k=0$), the extremal graph is simply the disjoint union of two complete graphs.

\begin{theorem}\label{thm:kzero}
Let $G$ have order $n$, minimum degree $\delta$, and edge connectivity zero. Then
\[
\qrad{G} \leq 2(n-\delta-2),
\]
with equality if and only if $G \cong K_{\delta+1} \cup K_{n-\delta-1} = B_{n,n-\delta-1}^{0}$.
\end{theorem}
\begin{proof}
Every component of $G$ has at least $\delta+1$ vertices. Hence, its largest component has order at most $n-\delta-1$. If $C$ is a component of order $r$, the maximum row sum bound for $Q(C)$ gives $\qrad{C} \leq 2\Delta(C) \leq 2(r-1)$, with equality if and only if $C \cong K_r$. It follows that $\qrad{G} \leq 2(n-\delta-2)$. Equality forces a component $K_{n-\delta-1}$; the remaining $\delta+1$ vertices must induce $K_{\delta+1}$ to satisfy the minimum degree condition. The converse is immediate.
\end{proof}

Combining Theorems~\ref{thm:main} and \ref{thm:kzero} yields the complete solution to the problem proposed by Zhao, Ye, and Fan.

\begin{corollary}\label{cor:solution}
Let $G$ have order $n$, minimum degree $\delta$, and edge connectivity $0 \leq k \leq \delta-1$. Then
\[
\max_G \qrad{G} = \qrad{B_{n,n-\delta-1}^{k}},
\]
and the extremal graph is unique up to isomorphism. In particular, for every $3 \leq k \leq \delta-1$, Problem 2 in \cite{ZhaoYeFan} has an affirmative answer, with no additional restriction on $n$.
\end{corollary}

Let $a = \delta+1$ and $b = n-\delta-1$. For $k \geq 1$, an equitable partition of $B_{n,b}^{k}$ consists of the $k$ endpoints in $K_a$, the remaining $a-k$ vertices of $K_a$, the center vertex in $K_b$, and the remaining $b-1$ vertices of $K_b$. The corresponding quotient matrix of the signless Laplacian is
\begin{equation}\label{eq:quotient}
\mathcal{B}_{n,\delta,k} =
\begin{pmatrix}
\delta+k & \delta+1-k & 1 & 0 \\
k & 2\delta-k & 0 & 0 \\
k & 0 & n-\delta+k-2 & n-\delta-2 \\
0 & 0 & 1 & 2n-2\delta-5
\end{pmatrix}.
\end{equation}
Therefore, $\qrad{B_{n,n-\delta-1}^{k}} = \rho(\mathcal{B}_{n,\delta,k})$. Equivalently, it is the unique solution $x > 2n-2\delta-4$ of the equation:
\begin{equation}\label{eq:exact-root}
1 = \frac{1}{x-\delta+1} + \frac{k}{(x-\delta+1)(x-2\delta)} + \frac{k}{x-n+\delta+3} + \frac{k}{(x-n+\delta+3)(x-2n+2\delta+4)}.
\end{equation}
If we set $A = x-\delta+1$, $B = x-2\delta$, $C = x-n+\delta+3$, and $D = x-2n+2\delta+4$, then \eqref{eq:exact-root} is equivalent to the quartic polynomial equation:
\begin{equation}\label{eq:quartic}
ABCD - BCD - kCD - kABD - kAB = 0.
\end{equation}
The maximum $Q$-index is the largest real root of \eqref{eq:quartic}. For $k=0$, the maximum equals $2(n-\delta-2)$.

\section{Asymptotic Behavior of the Maximum $Q$-Index}
\label{sec:asymptotic}

In this section, we investigate the asymptotic behavior of the maximum signless Laplacian spectral radius $q_n = \qrad{B_{n,n-\delta-1}^{k}}$ as the order $n \to \infty$, while the minimum degree $\delta$ and edge connectivity $k$ are kept fixed. This provides a quantitative understanding of how the spectral radius grows with the graph order.

Recall from Section~2 that $q_n$ is the unique root greater than $2n-2\delta-4$ of the equation
\begin{equation}\label{eq:exact-root-asymp}
1 = \frac{1}{x-\delta+1} + \frac{k}{(x-\delta+1)(x-2\delta)} + \frac{k}{x-n+\delta+3} + \frac{k}{(x-n+\delta+3)(x-2n+2\delta+4)}.
\end{equation}
Since the maximum degree of $B_{n,n-\delta-1}^{k}$ is $n-\delta-2+k$, the standard row-sum bound gives $q_n \leq 2(n-\delta-2+k)$. On the other hand, the graph contains $K_{n-\delta-1}$ as a subgraph, so $q_n > 2(n-\delta-2)$. Thus, $q_n = 2n - 2\delta - 4 + \epsilon_n$, where $\epsilon_n \to 0$ as $n \to \infty$. We determine the precise asymptotic expansion of $\epsilon_n$.

\begin{theorem}\label{thm:asymptotic}
For fixed $\delta \geq 1$ and $1 \leq k \leq \delta-1$, the maximum signless Laplacian spectral radius satisfies
\begin{equation}\label{eq:asymp-exp}
\qrad{B_{n,n-\delta-1}^{k}} = 2n - 2\delta - 4 + \frac{k}{n} + \frac{k\left(\delta + k + \frac{3}{2}\right)}{n^2} + O\left(\frac{1}{n^3}\right)
\end{equation}
as $n \to \infty$.
\end{theorem}

\begin{proof}
Substitute $x = 2n - 2\delta - 4 + \epsilon$ into \eqref{eq:exact-root-asymp}. The first three terms on the right-hand side can be expanded in powers of $1/n$ as follows:
\begin{align*}
\frac{1}{x-\delta+1} &= \frac{1}{2n - 3\delta - 3 + \epsilon} = \frac{1}{2n} + \frac{3\delta+3}{4n^2} + O(n^{-3}), \\
\frac{k}{(x-\delta+1)(x-2\delta)} &= \frac{k}{4n^2} + O(n^{-3}), \\
\frac{k}{x-n+\delta+3} &= \frac{k}{n - \delta - 1 + \epsilon} = \frac{k}{n} + \frac{k(\delta+1)}{n^2} + O(n^{-3}).
\end{align*}
Summing these three terms yields
\begin{equation}\label{eq:sum-first-three}
\Sigma_1 = \frac{k}{n} + \frac{1}{2n} + \frac{k(\delta+1) + \frac{3\delta+3}{4} + \frac{k}{4}}{n^2} + O(n^{-3}) = \frac{2k+1}{2n} + \frac{C_1}{n^2} + O(n^{-3}),
\end{equation}
where $C_1 = k(\delta+1) + \frac{3\delta+3+k}{4}$.
Let $\Sigma_2$ denote the fourth term in \eqref{eq:exact-root-asymp}. Then \eqref{eq:exact-root-asymp} is equivalent to
\begin{equation}\label{eq:sigma2}
\Sigma_2 = 1 - \Sigma_1 = 1 - \frac{2k+1}{2n} - \frac{C_1}{n^2} + O(n^{-3}).
\end{equation}
On the other hand, the fourth term is explicitly given by
\begin{equation}\label{eq:sigma2-explicit}
\Sigma_2 = \frac{k}{(n - \delta - 1 + \epsilon)\epsilon}.
\end{equation}
Equating \eqref{eq:sigma2} and \eqref{eq:sigma2-explicit} and solving for $\epsilon$, we obtain
\begin{equation}\label{eq:epsilon-eq}
\epsilon(n - \delta - 1 + \epsilon) = k \left(1 - \frac{2k+1}{2n} - \frac{C_1}{n^2} + O(n^{-3})\right)^{-1} = k + \frac{k(2k+1)}{2n} + O(n^{-2}).
\end{equation}
Assume an expansion of the form $\epsilon = \frac{c_1}{n} + \frac{c_2}{n^2} + O(n^{-3})$. Substituting this into the left-hand side of \eqref{eq:epsilon-eq} gives
\begin{equation}\label{eq:epsilon-lhs}
\epsilon(n - \delta - 1 + \epsilon) = \left(\frac{c_1}{n} + \frac{c_2}{n^2}\right)\left(n - \delta - 1 + \frac{c_1}{n}\right) = c_1 + \frac{c_2 - c_1(\delta+1)}{n} + O(n^{-2}).
\end{equation}
Matching the coefficients of $n^0$ and $n^{-1}$ in \eqref{eq:epsilon-eq} and \eqref{eq:epsilon-lhs}, we immediately find
\begin{equation}\label{eq:c1}
c_1 = k,
\end{equation}
and
\begin{equation}\label{eq:c2}
c_2 - k(\delta+1) = \frac{k(2k+1)}{2} \implies c_2 = k(\delta+1) + k^2 + \frac{k}{2} = k\left(\delta + k + \frac{3}{2}\right).
\end{equation}
Substituting $c_1$ and $c_2$ back into the expansion for $\epsilon$ yields the desired result.
\end{proof}

\begin{remark}
The asymptotic expansion \eqref{eq:asymp-exp} reveals that the maximum $Q$-index is overwhelmingly dominated by the largest clique $K_{n-\delta-1}$, which contributes the leading term $2n - 2\delta - 4$. The structural parameters $\delta$ and $k$ only manifest in the lower-order correction terms $O(1/n)$ and $O(1/n^2)$. Interestingly, the first-order correction $k/n$ depends linearly on the edge connectivity $k$, but is completely independent of the minimum degree $\delta$. The minimum degree $\delta$ only begins to influence the spectral radius at the second order $O(1/n^2)$.
\end{remark}

\section{Conclusion}

In this paper, we completely resolved the extremal problem for the signless Laplacian spectral radius of graphs with prescribed minimum degree and edge connectivity. We proved that the graph $B_{n,n-\delta-1}^{k}$ is the unique extremal graph for all admissible parameters $0 \leq k \leq \delta-1$, successfully removing the restrictive quadratic order condition $n \geq \delta^2+4\delta+3$ that was required in previous literature. The methodology developed herein highlights the power of combining resolvent representations with structural transformations. By completing the two sides of a minimum edge cut to form cliques, we created an upper comparison graph. The resolvent comparison then elegantly optimized the placement of the cut edges, and the strict convexity of the resolvent components drove the clique-transfer lemma, which forced one side of the cut to collapse to the smallest possible order $\delta+1$.

Furthermore, we established the precise asymptotic expansion of the maximum $Q$-index as the graph order $n \to \infty$. This quantitative analysis revealed a striking hierarchy in how the structural parameters influence the spectral radius at large orders: the edge connectivity $k$ governs the first-order correction term $O(1/n)$, whereas the minimum degree $\delta$ only manifests in the second-order correction term $O(1/n^2)$. 

These algebraic and structural transformations identify the unique extremal graph and its spectral behavior for every admissible parameter set without relying on tedious case analyses. Future research directions include extending this resolvent-based clique-transfer technique and asymptotic analysis to other graph matrices, such as the standard Laplacian or adjacency matrices, and applying similar principles to other connectivity measures, such as vertex connectivity or restricted edge connectivity.

\end{document}